\RequirePackage{fix-cm}
\documentclass[english,11pt]{article}  
\usepackage{graphicx}
\usepackage{amsmath,amsfonts,latexsym,amssymb,amscd, graphicx}
\usepackage{latexsym,enumerate,epsfig,listings}
\usepackage{amsopn,amstext}
\usepackage[latin1]{inputenc}
\usepackage{verbatim}
\usepackage{graphicx}
\usepackage{hyperref}
\usepackage[dvipsnames, usenames]{color}
\usepackage{anysize}

\numberwithin{equation}{section}
\newtheorem{theorem}{Theorem}[section]

\newtheorem{lemma}[theorem]{Lemma}

\newtheorem{remark}[theorem]{Remark}

\renewcommand{\P}{{\mathbb P}}
\newcommand{\E}{{\mathbb E}}
\newcommand{\R}{{\mathbb R}}

\newcommand{\ZZ}{{\mathbb Z}}
\newcommand{\NN}{{\mathbb N}}

\newcommand{\0}{{\mathbf 0}}

\newcommand{\cc}{{\mathbf c}}

\newcommand{\yy}{{\mathbf y}}

\newcommand{\xx}{{\mathbf x}}

\newcommand{\ee}{{\mathbf e}}

\newcommand{\Z}{{\mathbb Z}}

\begin{document}

\title{Attractiveness of Brownian Queues in Tandem}



\author{Eric A. Cator         \and
  Sergio I. L\'opez \and Leandro P. R. Pimentel 
}



\maketitle

\begin{abstract}
	
Consider a sequence of $n$ bi-infinite and stationary Brownian queues in tandem. Assume that the arrival process entering in the first queue is a zero mean ergodic process. We prove that the departure process from the $n$-th queue converges in distribution to a Brownian motion as $n$ goes to infinity. In particular this implies that the Brownian motion is an attractive invariant measure for the Brownian queueing operator. Our proof exploits the relationship between the Brownian queues in tandem and the last-passage Brownian percolation model, developing a coupling technique in the second setting. The result is also interpreted in the related context of Brownian particles acting under one sided reflection.\\

\textbf{Subclass} 60K25, 60K35.
\end{abstract}

\section{Introduction}\label{intro}

Tandem queues systems (TQ) are classical models in queueing theory consolidated from many decades of research and generalized to stochastic networks with diverse structures. A tandem queue is a system of queues where there is an initial arrival process $A^1$ and a sequence $\{ S^n \}_{n \geq 1}$ of service processes, all independent. The system is defined recursively: the initial queue is fed from the arrival process $A^1$, and has departures determined by the service process $S^1$. For $n \geq 2$, the arrival process for the $n$-th queue is defined as the departure process of the ($n$-$1$)-th queue and the departures are determined by the service process $S^n$. One fundamental result in queueing theory is Burke's theorem, which states that, given a Poisson process as arrival and an independent Poisson process as service (where the service intensity is strictly larger than the arrival one), the departure process is a Poisson process. This type of result, where there is an invariant law of the process under the queueing operator, is known as an \textit{Output theorem} in the literature, and it allows to compute explicitly many features of tandem queues systems. 

It is natural to consider the convergence of the departure process law from the $n$-th queue, as $n$ goes to infinity, when the initial arrival process is arbitrary. This was answered in \cite{MP} in the case when the service processes are Poisson: there is convergence to a Poisson process, under weak conditions on the initial arrival process. In  \cite{Pra} the result was generalized to the case when the service processes are not Poisson but independent and identically distributed. In this work we study the same question when the service processes are Brownian motions. 

Let us start by introducing the Brownian Tandem Queues (TQ). We follow the notation introduced in \cite{OY}. For real and continuous functions $f \in \cal C (\R)$, set $f(x,y):=f(y)-f(x)$.  Let $a=(a(x)\,,x\in\R)$ denote some  continuous arrival process and for $\mu >0$ define the service process by $s^{(1)}(x):=\mu x-B^{(1)}(x)$, where $B^{(1)}=( B^{(1)} (x) \, , x \in \R)$ is a two-sided Brownian motion independent of $a$. The queue length process is defined as 
\begin{equation}\label{def:queue}
 q^{(1)}(x):=\sup_{z \leq  x}    \left\{  a(z,x)-s^{(1)} (z,x)   \right\} .
\end{equation}
In order for $q^{(1)}$ to be stable (positive recurrent), we impose that the service process $s^{(1)}$ has a drift larger than that of the arrival process. We do this by requiring
\[ \lim_{x\to -\infty} \frac{a(x)}{x} = 0 \mbox{ and } \lim_{x\to \infty} \frac{a(x)}{x} = 0.\]
The departure process is defined by 
\begin{equation}\label{def:dep}
d^{(1)}(x,y):=a(x,y) - q^{(1)}(x,y) ,
\end{equation}
with the convention $d^{(1)}(0)=0$, and hence we put $d^{(1)}(x):= d^{(1)}(0,x)$.

The tandem queue model, in words, consists of a line of queues, where each queue uses as input (arrival) process the output (departure) process of the queue that is just in front of it in the line. In this context we have an initial arrival process $a$ and service processes $\{ s^{(n)}  \}_{n \in \NN} $ where $s^{(n)}(x)=\mu x-B^{(n)}(x)$ and $\left\{B^{(n)}\,:\,n\in\NN\right\}$ is a collection of independent (two-sided) Brownian motions. 
One can define inductively the queue length and the departure process of the $n$-th Brownian queue. Assume that the departure process $(d^{(n)} (x) : x \in \mathbb R)$ is already defined. Then we can define the queue length process of the $n+1$-th Brownian queue as
$$q^{(n+1)}(x):=\sup_{z\leq x}\left\{d^{(n)}(z,x)-\mu (x-z)+B^{(n+1)}(z,x)\right\}\,,$$
and the departure process from the $n$-th Brownian queue 
$$d^{(n+1)}(x,y):=d^{(n)}(x,y)- q^{(n+1)}(x,y)\, ,$$
with the similar convention $d^{(n+1)}(0)=0$ and $d^{(n+1)}(x):= d^{(n+1)}(0,x)$.

A measure on the space of continuous arrival functions with zero drift is called invariant for the queueing operator (in equilibrium), if the departure process has the same law as the arrival process. For the Brownian queue operator, the measure induced by an independent standard Brownian motion $B$ is an invariant (ergodic) measure \cite{OY}. Our result is the uniqueness of such a measure, by proving attractiveness:
\begin{theorem}\label{attract}
	Start the process of queues in tandem with a zero mean ergodic arrival process. Then 
	\begin{equation}\label{ergodic}
	\lim_{n\to\infty}d^{(n)}\stackrel{dist.}{=}B\,.
	\end{equation}
\end{theorem}
In our way to prove Theorem \ref{attract} we will only use that $B$ is an invariant ergodic measure for the queue system. Uniqueness will follow from our method. 

Essential for our proof is the connection of the Brownian TQ model to two related Brownian models, namely the Brownian Last Passage Percolation (LPP) and the Totally Asymmetric Brownian Exclusion Process (TABEP). We will introduce these models in Section \ref{sec:rel}, and point out the relationships between the three models.

All of these models have been previously studied, and the connection between them has been known for a while. Hambly, Martin and O'Connell \cite{HMO} defined the LPP Brownian model and derived concentration results for the associated Brownian growth model. The related Brownian particle system model has been studied in \cite{FSW1,FSW2}: particles are driven by Brownian motions and each particle is reflected (only) on its left closest particle. While models of Brownian motions interacting by exclusion on the real line have been an active research topic \cite{Ich,KPS,PaPi}, Ferrari, Spohn and Weiss successfully constructed a strong version of a two-sided system with an infinite amount of particles in a stationary regime \cite{FSW2}, governed by an asymmetric Skorokhod's type reflection, easily related to the LPP model. They accomplished it by some technique resembling Loynes' stability theorem for $G/G/1$ queues \cite{Loy}, and studied the finite-dimensional distributions of the system, characterized in terms of the Airy process. For simplicity, we name this Brownian particle system as the Totally Asymmetric Brownian Exclusion Process (TABEP), as suggested by P. A. Ferrari. The queueing model related to this particle process is exactly our Brownian TQ.

\subsection{Contribution}

In this article, we first revisit the connection between these three models: the LPP Brownian model, the TABEP process and the TQ Brownian system. The relationship between the LPP model and the TABEP process is mentioned in \cite{FSW2} while the relation between the LPP model and the TQ system is described in \cite{OY}. This is completely analogous to the known relationship between the standard Markovian Tandem Queues, LPP on ${\mathbb Z}^2$ with exponential weights and the TASEP (Totally Asymmetric Exclusion process). For the sake of completeness, these models are presented in Section \ref{sub:disc}.

Relying on these relations, we prove a result concerning the uniqueness of the invariant measure for the Brownian queueing operator, by proving attractiveness to that measure. In words, if we start with some zero mean ergodic process as initial arrival process and let it pass through the Brownian queues in tandem then the departure process from the $n$-th queue converges in distribution to a Brownian motion as $n$ goes to infinity. This is precisely stated in Theorem \ref{attract}.

For this purpose, we only use that the invariant measure under the queueing operator is known \cite{HW}: it is the random measure associated to the Brownian motion. The method of proof is a coupling technique developed in the LPP Brownian setting: starting with two different initial arrival processes (called \textit{mass profiles} in the LPP Brownian model) we use the same service processes (the random environment in the LPP context) to define the coupled evolution. Then we can prove that the difference between the associated departure processes (mass profiles) at each stage of the tandem is converging to zero on compact sets. This is our main result: Theorem \ref{lppattract}, which implies the desired conclusion in the queueing context, Theorem \ref{attract}. We point out that this result can also be translated to an attractiveness result for a semi-infinite TABEP system, see Theorem \ref{TABEPatt}.

A key step in the method involves local comparison techniques which allow us to bound the difference between mass profiles in terms of the so-called \textit{exit points} in the LPP literature. This implies that it is only necessary to control the exit points for a given system (done in Lemma \ref{control}) and then to control the difference between the exit points defined for each of the coupled systems. These exit points are naturally defined in the LPP context but we give an interpretation in the queueing setting in the following. First, consider an arbitrary initial arrival process and a single node Brownian queue. The exit point associated to time $x$ is the last time $Z(x,1)$ before time $x$ when the Brownian queue was empty. Given node $n$ of a tandem Brownian queue system and some time $x$, define $I_{n-1}(x,n)$ as the last time the $n$-th queue was empty before time $x$, then $I_{n-2}(x,n)$ to be the last time the ($n$-$1$)-th queue was empty before time $I_{n-1}(x,n)$, and so on, until we find the exit point $Z(x,n)=I_0(x,n)$. Hence, the exit time can be found from this iterative process of marking the beginning of the current excursion of the queue in each stage of the tandem system. This property is described in more detail in Subsection \ref{sub:exit}.
 
Our method of proof differs substantially from the methods developed for discrete valued queueing systems: In \cite{MP} a coupling between the departure times in every step of the tandem queue of each user is accomplished, while in \cite{Pra} the waiting times of each user in every node of the tandem queue system are considered for the coupling.

A rather simplified version of this result was presented in \cite{Lop} where, using a path coupling of the departures processes, a non-stationary and one-sided (in time) system is studied with some particular initial conditions. Those techniques are non applicable to the current bi-infinite stationary setting.

\subsection{Structure of the paper}

In Section $2$, we first review the classical discrete models. Then we define the Totally Asymmetric Brownian Exclusion Process (Subsection \ref{sub:TABEP}) and the Last Passage Percolation System (Subsection \ref{sub:LPP}). In each of these Subsections, our result is stated in the corresponding context (Theorems \ref{TABEPatt} and \ref{lppattract}) and the explicit relations between the models are shown. In Subsection \ref{sub:LPP} the coupled dynamics are defined. In Section 3, we first present the definition of exit points and the results concerning its control (Subsection \ref{sub:exit}) and then proceed to show the comparison results and the proof of Theorem \ref{lppattract} (Subsections \ref{sub:comparison} and \ref{sub:proof}).

\section{The discrete models}\label{sub:disc}

In this Section we review some fundamental relationship between the classical Markovian Tandem queue model (TQ) with the exponential last-passage percolation model (LPP) and the Totally Asymmetric Simple Exclusion Process (TASEP).\\ 

Assume that we have $K$ Markovian queues in tandem working under a FIFO discipline. At time zero, the first queue starts working with $N$ users in the line while all the other queues are empty. Define a collection of rate one independent exponential random variables $\{ X(n,k) \}_{ n=1,...,N, \\ k=1,...,K}$ where $X(n,k)$ represents the service time of the $n$-th user at server $k$. Define $D(n,k)$ as the time where the $n$-th user exits the $k$-th server. Note that server $k$ only starts to serve user $n$ after user $n-1$ has exited server $k$ and the service from server $k-1$ to user $n$ has been finished. Then we have the following recurrence structure:
\begin{equation}\label{Lindley}
D(n,k)= X(n,k) + \max (  D(n,k-1), D(n-1,k) ),
\end{equation}
with boundary conditions $D(0,0)=0$ and $D(n,k)=0$ if $n<0$ or $k<0$. We will show how this structure is related with the aforementioned models.

Consider a collection of i.i.d. random variables $\{W_\xx\,:\,\xx\in (\ZZ_+)^2\}$ (also called weights), distributed according to an exponential distribution function of parameter one. In last-passage site percolation (LPP) models, each number $W_\xx$ is interpreted as the percolation (passage) time through  vertex $\xx=(x(1),x(2))$. For a  lattice vertex $\xx= (n,k)$ in $(\Z_+)^2$, denote $\Gamma(\xx)$ the set of all up-right oriented paths $\gamma=(\xx_0,\xx_1\dots,\xx_k)$ from $\0$ to $\xx$, i.e. $\xx_0=\0$, $\xx_k=\xx$ and $\xx_{j+1}-\xx_j\in\{\ee_1,\ee_2\}$, for $j=0,\dots,k-1$, where $\ee_1=(1,0)$ and $\ee_2=(0,1)$. The weight (or passage time) along $\gamma$ is defined as
$$W(\gamma):=\sum_{j=0}^{k} W_{\xx_i}\,.$$
The \emph{last-passage time} between $\0$ and $\xx$ is defined as
\begin{equation*}
L(\xx) \equiv L(n,k):=\max_{\gamma\in\Gamma(\xx)} W(\gamma)\,.
\end{equation*}
By the up-right path structure and the dynamic programming principle, we have the following Bellman equation:
\begin{equation}\label{Bellman}
L(n,k)= W(n,k) + \max (  L(n,k-1), L(n-1,k) ).
\end{equation}
This equation is the same as \eqref{Lindley} with the same boundary conditions, so the last passage percolation function is an equivalent way to describe the departure times from a tandem queue system.

Let us define the related interacting particle system. Let $\Omega$ be the space of binary sequences $\eta: \ZZ \rightarrow \{0, 1\}$. The elements $\eta$ in $\Omega$ will be configurations of particles. We will say that a configuration $\eta$ such that $\eta(x)=1$ has a particle at position $x$. If $\eta(x)=0$ we say that position $x$ is empty or that we have a \textit{hole} in that position. The dynamics are defined by the infinitesimal generator
\begin{equation*}\label{TASEP}
\mathcal L [f] (\eta )= \sum_{ x \in \mathbb Z}
\eta(x) (1 - \eta(x-1)) ( f ( \eta^{x,x-1}) - f(\eta)),
\end{equation*}         
where $\eta^{x,x-1}$ is defined as the configuration that is identical to $\eta$ except for the positions $x$ and $x-1$, where the original values are exchanged. The interpretation is the following: from each possible site $x$ we have a constant rate of jump of the particles (if there is no particle at site $x$, nothing happens). Once the clock at position $x$ rings, the particle in that place tries to jump to the site $x-1$ and this is accomplished if the site $x-1$ is empty, otherwise the jump is disregarded. This last condition emulates an exclusion principle, which is the reason that this process is known as the Totally Asymmetric Simple Exclusion Process. It is a standard microscopical model for transport, see for example \cite{Ferr}.



Finally, we show the relationship between the TASEP process and the Tandem queue model defined by \eqref{Lindley}. Let $\{ \eta_t : t \geq 0 \}$ be a TASEP process with initial configuration $\eta_0$. Assume the initial configuration $\eta_0$ is such that $\eta_0(x)= 1_{ [0, \infty)  } (x)$ for every $x \in \mathbb Z$, which means that all the particles are at the right of the origin in consecutive positions. Label each particle with its initial position, and define $x_{l}(t)$ to be the position of the $l$-th particle at time $t$ (so $x_l(0) = l$ for every $l \in \mathbb N$). Define
\begin{equation}\label{TASEPtoTQ}
q_l (t) := x_l(t)- x_{l+1}(t)-1,
\end{equation}
that is, the number of users in server $l$ at time is equal to the number of holes between particles $l$ and $l+1$ at time $t$. Note that \eqref{TASEPtoTQ} translates exactly the movement of particles in the exclusion process to the tandem queues dynamics: every time that the particle $l$ moves to the left, one user is entering into the $l^{\rm th}$ queue. Moreover, if the particle $l$ which is moving is not the first one, the number of users in the ($l$-$1$)-th queue diminish by one, so the user is leaving that queue. The exclusion property for particles translates into the restriction of having a non-negative number of users in each queue.
Consider now that holes are labeled in the starting configuration $\eta_0$: the hole at position $l <0$ will have label $-l$. The model is symmetric in particles and holes: one can think of holes traveling to the right which satisfy the exclusion property between them. Therefore, using \eqref{TASEPtoTQ}, we have another interpretation of the departure time $D(n,k)$: it is exactly the time when particle $k$ is exchanging position with hole $n$.\\

The previous presented relationships are known and studied, see \cite{Martin}. In the last two decades great progress has been made for the LPP models and this has given insight to an important question originally posed in queueing theory: the asymptotic distribution of the departure time of the $n$-th user in line from the $m$-th queue (its order in the line of queues), when the whole system starts empty, by making $m$ and $n$ grow to infinity while keeping fixed the ratio between them \cite{Gly,Sep}. On the other hand, strong results from queueing theory concerning the existence and attractiveness of invariant measures under the queueing operator \cite{Mai,Pra} have been used to shed light on difficult questions concerning LPP models, as for example the existence of semi-infinite geodesics and Busemann functions for the lattice model in $\mathbb{Z}^2$ with general distributed weights, see \cite{GRS2,GRS1}.

\section{Convergence of the Brownian Models}\label{sec:rel}

Theorem \ref{attract} states our convergence result for the Brownian TQ. In this Section and the next one we will restate basically the same result in the context of two different Brownian models.

\subsection{Convergence in the Totally Asymmetric Brownian Exclusion Process}\label{sub:TABEP}

Consider a semi-infinite system of Brownian interacting particles defined for all real times $x$. Take some stationary, ergodic and continuous process $\{X^{(0)}(x): x \in \mathbb{R} \}$ and define $X^{(0)}(x)$ as the position of the leftmost particle at time $x$. We introduce a collection $\{B^{(n)}:n\geq 1\}$ of independent two-sided standard Brownian motions. Then, for $n \geq 1$, define
\begin{equation}\label{TABEPsta}
X^{(n)}(x) = \sup_{ y \leq x} (X^{(n-1)}(y)  + B^{(n)}(x) - B^{(n)}(y)), \qquad x \in \mathbb R \, .
\end{equation}
The system $\{X^{(n)}(x): x \in \mathbb{R} \}_{ n \geq 0}$ will be called the Totally Asymmetric Brownian Exclusion Process (TABEP) with leftmost particle $X^{(0)}$. By definition, the order of the particles is preserved: $X^{(0)}(x) \leq X^{(1)}(x) \leq ... $ for every real time $x$ (choosing $y=x$ in the argument of the supremum in \eqref{TABEPsta} shows that $X^{(n-1)}(x) \leq X^{(n)}(x)$). Note that \eqref{TABEPsta} implies that the TABEP is Markovian in $n$: conditionally on the information of the process $X^{(n)}$, the process $X^{(n+1)}$ is independent from the collection $\{X^{(k)} \}_{ k=1,...,n-1}$. These two properties can be combined to give an informal interpretation: the $n$-th particle 
is obtained by reflecting an independent Brownian motion to its left-side neighbor (the ($n$-$1$)-th particle) and this is the only possible interaction between particles (note that a particle does not notice the particles to the right of it).

 A sufficient condition to have a well-defined system is that for some positive constant $\mu$, $X^{(0)}$ satisfies
\begin{equation}\label{leftmost}
\liminf_{x\to-\infty}\frac{ X^{(0)} (x)}{x}\geq \mu\,\,\,\mbox{ and that }\,\,\,\limsup_{x\to\infty}\frac{ X^{(0)} (x)}{x}\leq \mu\,. 
\end{equation}
Note that the whole system is time stationary: one can prove inductively that the distribution of $X^{(n)}(x)$ does not depend on $x$, for every $n\geq 0$. 

Let us remark that the system defined above is a two-sided time stationary extension of a TABEP system with initial positions, defined by Ferrari, Spohn and Weiss \cite{FSW2}. In that work, they considered the particular case of initial positions where the starting positions of the particles are given by a rate $\mu$ Poisson process on $[0,\infty)$ and the left-most particle is given by
 $$X^{(0)}(x) = B^{(0)}(x) + \mu x,$$
 where $\{ B^{(0)}(x): x \in \mathbb R \}$ is a Brownian motion. Using Burke's theorem for Brownian motion \cite{OY}, they constructed a stationary bi-infinite system of ordered particles
$$ ... \leq  X^{(-1)} (x) \leq X^{(0)} (x) \leq X^{(1)}(x) \leq ...     \qquad  \forall x \geq 0$$
where each particle has the distribution of a standard Brownian motion (where the initial position is not zero) and, for each positive time $x$, the set of positions is distributed as a rate $\mu$ Poisson process on the line. 

Now we show the relation with the tandem Brownian queues. Consider a TABEP system $\{ X^{(n)} \}_{n \geq 0}$, defined by \eqref{TABEPsta}. Define the arrival process $a(x):=\mu x  -X^{(0)}(x)$ and the service processes $s^{(n)}(x):=\mu x -B^{(n)}(x) $ for each $n \geq 1$ (note that $a$ has zero drift and $s^{(n)}(x)$ has positive drift $\mu$). Then the associated first queue length process is given by
$$ q^{(1)}(x) = \sup_{ y \leq x}  (X^{(0)}(y) - X^{(0)}(x)+ B^{(1)}(x) - B^{(1)}(y) ) \quad \forall x \in \mathbb{R},$$
the first departure process is 
$$  d^{(1)}(x) = q^{(1)}(0) +X^{(0)}(0) + \mu x - \sup_{ y \leq x } (X^{(0)}(y)  + B^{(1)}(x) - B^{(1)}(y) )    \quad \forall x \in \mathbb{R}, $$
and, by \eqref{TABEPsta}, we conclude that $d^{(1)}(x) = X^{(0)}(0) + q^{(1)}(0) + \mu x  - X^{(1)}(x)$ (we are using the convention $d^{(1)}(0)=0$).
  
Analogous formulae hold for any $n\geq 1$, by induction: Suppose now that for a fixed natural $k$ we have
$$ d^{(k)} (x)= X^{(0)}(0) + \sum_{i=1}^{k} q^{(i)}(0)  + \mu x  - X^{(k)}(x),   \quad \forall x \in \mathbb R .$$
Then
$$ q^{(k+1)}(x) = \sup_{ y \leq x} (  d^{(k)} (y,x) - s^{(n)}(y,x) ) = \sup_{ y \leq x}  (X^{(k)}(y) - X^{(k)}(x)+ B^{(k)}(x) - B^{(k)}(y) ),  \quad \forall x \in \mathbb{R}. $$
Since $d^{(k+1)}(0)=d^{(k)}(0)=0$, this implies that
\begin{eqnarray}
 d^{(k+1)} (x) &=& d^{(k)} (x) - q^{(k+1)}(x) + q^{(k+1)}(0) \\
 &=& X^{(0)}(0) + \sum_{i=1}^{k+1} q^{(i)}(0)  + \sup_{ y \leq x}  (X^{(k)}(y) + B^{(k)}(x) - B^{(k)}(y) )    \quad \forall x \in \mathbb{R},
\end{eqnarray}
where we also used the induction hyphotesis. By \eqref{TABEPsta} it follows that
$$ d^{(k+1)} (x)= X^{(0)}(0) + \sum_{i=1}^{k+1} q^{(i)}(0)  + \mu x  - X^{(k+1)}(x),   \quad \forall x \in \mathbb{R}.$$

An important remark is that, by using \eqref{TABEPsta}, we get that
\[ q^{(n)}(x) = X^{(n)}(x) - X^{(n-1)}(x),\]
so the distance between the particles  $n-1$ and $n$ is equal to the $n$-th queue length process at time $x$.
Thus \eqref{ergodic} is equivalent to Theorem \ref{TABEPatt} below.
\begin{theorem}\label{TABEPatt}
	Start a two sided TABEP with an ergodic process as the leftmost particle
	which satisfies \eqref{leftmost} for some positive constant $\mu$. Then the limit of the (centered) $n$-th particle converges to a two-sided Brownian Motion with drift $\mu$, that is
	\begin{equation}\label{TABEPlim}
	\lim_{n\to\infty}\left( X^{(n)}(x) -X^{(n)}(0) \right) \stackrel{dist.}{=} B(x)  + \mu x  \, .
	\end{equation}
\end{theorem}

\subsection{Convergence of the Brownian Last-Passage Percolation System}\label{sub:LPP}

In this section we define the elements of the theory of last-passage percolation systems \cite{CaPi12} with Brownian passage times, as developed in \cite{HMO}, and show its relationship with tandem Brownian queues. Let $\omega:=\left\{B^{(n)}\,:\,n\in\ZZ\right\}$ be a collection of i.i.d. two-sided Brownian motions. Define the order ``$<$'' in $\R\times\ZZ$ as the coordinate-wise order. For $\xx=(x,k)< \yy=(y,l)\in\R\times\ZZ$ denote $\Gamma(\xx,\yy)$ 
the set of all real increasing sequences $\gamma=(x=z_0\leq z_{1}\leq \dots \leq z_{l-k+1}=y)$. The passage time of $\gamma$ is defined as 
\begin{equation*}\label{passage} 
L(\gamma):=\sum_{i=0}^{l-k} B^{(k+i)}(z_i,z_{i+1})\,.
\end{equation*}
The last-passage time between $\xx$ and $\yy$ is given by
\begin{equation}\label{lastpassage} 
L(\xx,\yy):=\sup_{\gamma\in\Gamma(\xx,\yy)}L(\gamma)\,.
\end{equation}

The passage time of a path $\gamma$ can be seen as a continuous real valued process 
$X=(X(z): z \in \Gamma)$ where
$$ \Gamma = \{ z=(z_1,...,z_{l-k}) : x \leq z_1 \leq \dots \leq z_{l-k} \leq y \}  \subseteq \R^{l-k}.$$ 
Since $\Gamma$ is compact, by continuity, we have that the maximum is attained at some location. In \cite{LoPi} is proven that, for $\xx$ and $\yy$ fixed, the maximum is attained at a unique location with probability one. However, it is not true that this uniqueness holds simultaneously for all points $\xx, \yy \in \mathbb{R} \times \mathbb{N}$. 
To see an example, for  $x>0$ define  
$$ \mathcal{Z}(x) = \{ z \in [0,x] \, : \, B^{(0)}{(0,z)} + B^{(1)}{(z,x)} = L(\0, (x,1) ) \, \},$$
where $\0=(0,0)$. Put $W_x:= B^{(0)}(x) - B^{(1)}(x)$ and note that $z \in \mathcal{Z}(x)$ is equivalent to $W_{z} = \sup_{ u \in [0,x] } W_u$. Thus, by Levy's theorem, we have that
$$ \{ x \geq 0 \, : \, \# \mathcal{Z}(x)> 1 \}  \stackrel{dist.}{=} 
\{ x \geq 0 \, : \, \sqrt{2} \, l_x \textrm{ is strictly increasing} \}, $$
where $l_x$ is the local time of a standard Brownian motion. \\




We will call the \textit{geodesic} (or the maximizer) between $\xx$ and $\yy$ to be the path $\gamma(\xx,\yy)$  such that 
$$L(\gamma(\xx,\yy))=L(\xx,\yy)\, .$$ 

To introduce the last-passage percolation system we consider an initial profile $\nu=(\nu(x)\,,x\in\R)$ such that $\nu(0)=0$ and
\begin{equation}\label{asymp}
\liminf_{ y \rightarrow  -\infty} \frac{ \nu (y) }{y} >0 ,
\end{equation}
and define the (discrete time) evolution of $\nu$ as the Markov process $(M^{(n)}_\nu\,,\,n\geq 0)$, where $M^{(0)}_{\nu}=\nu$,
\begin{equation}\label{system} 
L_\nu(x,n):=\sup_{z\leq x}\left\{\nu(z)+L\left((z,1),(x,n)\right)\right\}\,\,\mbox{ and }\,\,M^{(n)}_\nu(x):=L_\nu(x,n)-L_\nu(0,n)\,.
\end{equation}
The Markov property follows from the following fact: for all $n\geq 1$ and $k\in\{0,\dots,n-1\}$
\begin{equation}\label{markov}
L_\nu(x,n) - L_\nu(0,k) =\sup_{z\leq x}\left\{M^{(k)}_\nu(z)+L\left((z,k+1),(x,n)\right)\right\}\,,
\end{equation}
which is an application of the dynamic programming principle. This is a graphical construction of the process where the space-time random environment is given by the collection of Brownian motions $\omega= \left\{  B^{(n)}\,:\,n\in\ZZ \right\}$. The variational formula expresses the profile at time $n$ as a function of the profile at time $k<n$ plus some strip of the space-time environment which is independent of the profile at time $k$. We note that this construction allows us to run the last-passage percolation system, started with two arbitrary initial profiles $\nu_1$ and $\nu_2$, simultaneously with the same environment $\omega$ (basic coupling). Formally speaking, we define the joint process $(M^{(n)}_{\nu_1},M^{(n)}_{\nu_2})_{n\geq 0}$ by setting
\begin{equation}\label{basic}
(x,n)\,\mapsto\,\left\{\begin{array}{ll} L_{\nu_1}(x,n):=\sup_{z\leq x}\left\{\nu_1(z)+L\left((z,1),(x,n)\right)\right\}\,,&\\
L_{\nu_2}(x,n):=\sup_{z\leq x}\left\{\nu_2(z)+L\left((z,1),(x,n)\right)\right\}\, ,&
\end{array}
\right.
\end{equation}
and putting $M^{(n)}_{\nu_i}(x):=  L_{\nu_i}(x,n)-  L_{\nu_i}(0,n)$ for $x$ real and $i=1,2$. Notice that $L\left((z,1),(x,n)\right)$ is a function that only depends on $\omega$.\\

The analogy with the queue system is as follows. Assume that $\nu(x)$ has drift $\mu$ and take 
\begin{equation}\label{LPPtTQ1}
a(x)=\mu x - \nu(x)\,\,\mbox{ and }\,\,s^{(n)}(x):=\mu x-B^{(n)}(x)\,.
\end{equation}
Then
$$q^{(1)}(x):= \sup_{z\leq x}\left\{a(z,x)-s^{(1)}(z,x)\right\}=L_{\nu}(x,1)-\nu(x)$$
and
$$d^{(1)}(x):=a(x) + q^{(1)}(0) - q^{(1)}(x)=\mu  x - M_\nu^{(1)}(x)\,.$$  
From this, using definitions \eqref{def:queue}, \eqref{def:dep}, \eqref{system} and induction, one can check the analogous relation for all $n\geq 2$: 
$$q^{(n)}(x)= \sup_{z\leq x}\left\{  d^{(n-1)}(z,x)-s^{(n)}(z,x)\right\}=L_{\nu}(x,n)-L_{\nu}(x,n-1)$$
and
$$d^{(n)}(x)=a(x) + q^{(n)}(0) - q^{(n)}(x)=\mu  x - M_\nu^{(n)}(x)\,.$$  

Thus, \eqref{ergodic} and \eqref{TABEPlim} are consequences of \eqref{ergodiclpp} below. Define
$$B_\mu(x)=\mu x+B(x)\,,$$
where $B$ is a standard Brownian motion. Using the invariance of the Brownian measure under the queueing operator, it is immediate that $B_\mu$ is invariant:
\begin{equation*}
M_{\mu}^{(n)}\equiv M_{B_\mu}^{(n)}\stackrel{dist.}{:=}B_\mu\,,\mbox{ for all }n\geq 0\,.
\end{equation*}

The main contribution of this article is the next theorem, from which  \eqref{ergodic} (and \eqref{TABEPlim}) will follow.  
\begin{theorem}\label{lppattract}
	Let $\mu\in(0,\infty)$ and assume that, almost surely,
	\begin{equation}\label{exitcontrol}
	\liminf_{x\to-\infty}\frac{\nu(x)}{x}\geq \mu\,\,\,\mbox{ and }\,\,\,\limsup_{x\to\infty}\frac{\nu(x)}{x}\leq \mu\,.
	\end{equation}
	Consider the basic coupling $(M_{\nu}^{(n)},M^{(n)}_\mu)_{n\geq 0}$ constructed by running the last-passage percolation system, started with $\nu$ and $B_\mu$, simultaneously with the same environment $\omega=\left\{B^{(n)}\,:\,n\in\ZZ\right\}$. Then, for all compact $K\subseteq\R$ and $\epsilon>0$, 
	\begin{equation}\label{ergodiclpp}
	\lim_{n\to\infty}\P\left( \sup_{x\in K}|M_{\nu}^{(n)}(n,\mu^{-2}n+x)-M_\mu^{(n)}(n,\mu^{-2}n+x)|>\epsilon \right)=0\,.
	\end{equation}
\end{theorem}

It should be clear that an ergodic initial profile satisfies \eqref{exitcontrol} almost surely (note that in that case, we have translation invariance of the law of $M_\nu^{(n)}$ and $M_\mu^{(n)}$, so that we can get rid of the translation by $\mu^{-2}n$). We note that  \eqref{ergodiclpp} implies local convergence for initial profiles beyond the ergodic condition: one could take a deterministic profile satisfying \eqref{exitcontrol}. 

\section{Proofs}

\subsection{Shape Theorem and Exit Points}\label{sub:exit}
First proven in \cite{Bar,Gra}, using that $L( \0, (n,n))$ has the same law as the largest eigenvalue of a $n \times n$ GUE random matrix, the shape theorem below is presented by Hambly et al. \cite{HMO} as a consequence of concentration results for the Brownian directed percolation paths: 
\begin{equation}\label{shape} 
\lim_{n\to\infty}\frac{1}{n} L( (\0, (xn,tn)) \stackrel{a.s.}{=} 2 \sqrt{xt}.
\end{equation}
Note that, by Brownian scaling, 
\begin{equation}\label{scal}
\left\{L((\0, (rn,n))\,:\,r\in[0,x]\right\}\stackrel{dist.}{=}\left\{\sqrt{x}L( (\0,(sn,n))\,:\,s\in[0,1]\right\}\,.
\end{equation}

\begin{remark}\label{ctrpath}
	By Lemma $7$ in \cite{HMO}, there exist constants $c_1,c_2 \geq 0$ such that 
	\begin{equation}\label{auxconc}
	\mathbb{P} \Big( \Big| \frac{ L ( \0, (n,n))}{n} - 2  \Big| \geq 2 y \Big) \leq c_1 \exp\{-c_2 \, n (y- \epsilon_n)^2 \},
	\end{equation}
	for all $n \geq 0$, and $y > \epsilon _n$, where
	$$\epsilon_n := 2 - \frac{\E L (\0,(n,n))}{n} + \frac{1}{n^{1/4}}\,.$$
	Since $\epsilon_n \to 0$, we can choose $n$ large such that $\epsilon_n < 4^{-1} \delta$ and take $y=2^{-1} \delta$. This implies that there exist constants $c_3,c_4>0$ such that for all $\delta>0$ there exists $N>0$ such that  
	$$ \mathbb{P} \Big( \Big| \frac{ L ( \0, (n,n))}{n} - 2  \Big| \geq \delta \Big) \leq c_3 \exp\{-c_3 \, n \delta^2 \} ,$$
	for all $n \geq N$. We notice that a better upper bound could be produced by using the coupling method \cite{CaGr} to prove that
	$$\E|L( \0, (n,n))-2n|=O(n^{1/3})\,,$$
	which would imply that $\epsilon_n=O(n^{-1/4})$. For the Brownian last-passage percolation model we have all the ingredients necessary for the coupling method: we know explicitly the invariant regime and the shape function. 
\end{remark}

From now on we will treat $\nu$ as a fixed deterministic profile satisfying \eqref{exitcontrol}. Define the exit point from $(x,n)$ as
\begin{equation}\label{def:exit}
Z_{\nu} (x,n) = \sup \left\{  z \leq x : L_{\nu}  (x,n) =  \nu(z) + L( (z,1),(x,n) )  \right\}.
\end{equation}
We note that it is well defined. First, since we have the same asymptotic hypothesis \eqref{asymp}
on the profile $\nu$, one can use similar arguments as in Proposition $4.1$ of \cite{CaPi12} to prove that the function $L_{\nu}  (x,n)$ is well defined. By Brownian continuity, the map $z \rightarrow L( (z,1),(x,n) )$ is continuous, just as the profile $\nu$ (by hypothesis). Then the set
$ \left\{  z \in \mathcal C : L_{\nu}  (x,n) =  \nu(z) + L( (z,1),(x,n) )  \right\} $ is non empty for any compact set $C$. To prove that the supremum over $z \leq y$ can be restricted to some compact
set one can mimic the proof of Lemma $4.3$ in \cite{CaPi12}.\\

The name \textit{exit point} comes from the next geometrical interpretation in last-passage percolation: $Z_{\nu} (x,n)$ is the time before $x$ when the path which maximizes the quantity $\nu(z) + L( (z,1),(x,n) )$ leaves the initial profile $\nu$ (that can be visualized on the line $\{ (x,0) : x \in \mathbb R \}$) to percolate to the point $(x,n)$. 

The exit point \eqref{def:exit} can also be described in terms of the tandem queueing system. First, let us examine the interpretation for $Z_{\nu} (x,1)$. Let $z^*$ be in $\left\{  z \leq x : L_{\nu}  (x,1) =  \nu(z) + L( (z,1),(x,1) )  \right\}$. Then
$$ \nu(z^*) + L( (z^*,1),(x,1) \geq \nu(z) + L( (z,1),(x,1)  \quad \forall z \leq x ,$$
and, by \eqref{LPPtTQ1}, this implies that
$$  a(z)-s^{(1)}(z)  \geq  a(z^*)-s^{(1)}(z^*)  \quad \forall z \leq x .$$
In other words, 
$$ a(z^*,x) - s^{(1)}(z^*,x) \geq  a(z,x) - s^{(1)}(z,x) \quad \forall z \leq x ,$$
so $q^{(1)}(x)= a(z^*,x) - s^{(1)}(z^*,x)$ (by the definition \eqref{def:queue}). This implies that $q^{(1)} (z^*)=0$, so the value $Z_{\nu} (x,1)$ is the last time when the queue-length process $q^{(1)}$ was empty before time $x$. For $n$ arbitrary, using the expression \eqref{markov}, one can check that the value $Z_{\nu} (x,n)$ can be obtained inductively: 
let $I_{n-1}(x,n)$ be the last time when $q^{(n)}$ was empty before time $x$, then $I_{n-2}(x,n)$ is the last time when $q^{(n-1)}$ was empty before time $I_{n-1}(x,n)$, and so on, till we find the exit point $Z_{\nu}(x,n)=I_0(x,n)$. \\

In the next result we show that, in probability, the exit point is asymptotically sublinear.
\begin{lemma}\label{exit}
	Let $\mu\in(0,\infty)$ and assume \eqref{exitcontrol}. Then, for all $C\in \R$ and $\epsilon >0$,
	$$\lim_{n\to\infty}\mathbb{P} \left(n^{-1}|Z_\nu(\mu^{-2}n+C,n)|>\epsilon \right)=0\,.$$
\end{lemma}

\textbf{Proof:}\\

	By Brownian scaling \eqref{scal}, one can restrict the attention to $\mu=1$. For fixed $\delta>0$, take $B_{1+\delta}$ and construct $L_{1+\delta}$ and $L_\nu$ simultaneously using the basic coupling \eqref{basic}. Since  
	$$L\left((z,1),(n+C,n)\right)\leq L_{1+\delta}(n+C,n) - B_{1+\delta}(z)\,$$
	and 
	$$L\left((1,0),(n+C,n)\right)=L\left((1,0),(n+C,n)\right)+\nu(0)\leq L_\nu(n+C,n)\,$$
	(recall that $\nu(0)=0$), we have that 
	$$\left\{Z_\nu(n+C,n)\geq u \right\}=\left\{\exists\,z\in[u,n+C]\,:\,\nu(z)+L\left((z,1),(n+C,n)\right)=L_\nu(n+C,n)\right\}\,$$
	is contained in the event 
	$$\left\{\exists\,z\in[u,n+C]\,:\,B_{1+\delta}(z)-\nu(z)\leq L_{1+\delta}(n+C,n)-L\left((1,0),(n+C,n)\right)\right\}\,.$$
	By \eqref{exitcontrol} there exists $K_0>0$ such that $\nu(z)\leq (1+\delta/2)z$ for all $z>K_0$. Hence, if $u>K_0$ then $\left\{Z_\nu(n+C,n)\geq u \right\}$ is contained in the event 
	\begin{equation}\label{exit11}
	\left\{\exists\,z\in [u,n+C]\,:\,B_{2^{-1}\delta}(z)\leq L_{1+\delta}(n+C,n)-L\left((1,0),(n+C,n)\right)\right\}\,. 
	\end{equation}
	
	Now we recenter the Brownian motion with drift at position $u$ by writing 
	$$B_{2^{-1}\delta}(z):=B_{2^{-1}\delta}(u)+\bar B_{2^{-1}\delta}(z)\,,$$
	where $\bar B_{2^{-1}\delta}(z):= B_{2^{-1}\delta}(z)-B_{2^{-1}\delta}(u)$ for $z\geq u$. Notice that $\{ \bar B_{2^{-1}\delta}(z) : z \geq u  \}$ has the same distribution as the process $ \{ B_{2^{-1}\delta}(z) : z\geq 0 \}$ and it is independent of $B_{2^{-1}\delta}(u)$. Let 
	$$A(u):=B_{2^{-1}\delta}(u)+\min_{z\geq u} \bar B_{2^{-1}\delta}(z)\,.$$ 
	This minimum is well defined because $\bar B_{2^{-1}\delta}$ has a positive drift, and its distribution is given by minus an exponential random variable of parameter $2^{-1}\delta$ (its value will not play an important role when $n$ grows to infinity, since $\delta$ is fixed). Thus, by \eqref{exit11},
	\begin{equation}\label{exit12}
	\left\{Z_\nu(n+C,n)\geq u \right\}\subseteq \left\{A(u)\leq L_{1+\delta}(n+C,n)-L\left((1,0),(n+C,n)\right)\right\}   \qquad \forall u \leq n+C  \,. 
	\end{equation}
	The strategy is to show that if $u=\epsilon n$ we can choose $\delta>0$ such that the event on the r.h.s. of \eqref{exit12} has small probability. For $\epsilon_1>0$, to be defined later, we have that the event on the r.h.s. of \eqref{exit12} has probability bounded by     
	$$\P\Big(L (\left(1,0),(n+C,n)\right)-2n\leq -\epsilon_1 n\Big)+\P\Big(A(u)\leq L_{1+\delta}(n+C,n)-2n+\epsilon_1 n\Big)\,.$$
	By the shape theorem,
	$$\lim_{n\to\infty}\P\Big(L\left((1,0),(n+C,n)\right)-2n\leq -\epsilon_1 n\Big)=0\,.$$ 
	On the other hand, 
	\begin{eqnarray}\label{sum}
	\nonumber\P\Big(A(u)\leq L_{1+\delta}(n+C,n)-2n+\epsilon_1 n\Big)&\leq &\P\Big(2^{-1}\delta u-2\epsilon_1 n\leq L_{1+\delta}(n+C,n)-2n\Big)\\
	\label{exit13}&+&\P\Big(A(u)\leq 2^{-1}\delta u-\epsilon_1 n\Big)\,. 
	\end{eqnarray}
	We now use the result in Section 4 of \cite{OY}, where it is shown (in our notation) that $L_\lambda(0,n) - L_\lambda(0,0)$ (this is the vertical increment) is distributed as the sum of $n$ \emph{independent} exponential random  variables, each with expectation $1/\lambda$. We already know that $x\mapsto L_\lambda(x,n)-L_\lambda(0,n)$ (the horizontal increment) is distributed as Brownian motion with drift $\lambda$. This shows us how to recenter $L_{1+\delta}(n+C,n)$:
	$$\P\Big(2^{-1}\delta u-2\epsilon_1 n\leq L_{1+\delta}(n+C,n)-2n\Big)=\P\left(\Delta -2\epsilon_1 n\leq L_{1+\delta}(n+C,n)-\left((1+\delta)+\frac{1}{1+\delta}\right)n\right)\,,$$
	where 
	\begin{eqnarray}
	\nonumber\Delta &:=& 2n-\left((1+\delta)+\frac{1}{1+\delta}\right)n+\frac{\delta}{2}u\\
	\nonumber &=&-\frac{\delta^2}{(1+\delta)}n+\frac{\delta}{2}u\\
	\nonumber &>& \left(\frac{\delta}{2}\frac{u}{n}-\delta^2\right)n\,.
	\end{eqnarray}
	If $u=\epsilon n$ and we pick $\delta:=4^{-1}\epsilon$, we get the next lower bound for $\Delta$:
	\begin{eqnarray}
	\nonumber\Delta &>& \left(\frac{\delta}{2}\epsilon-\delta^2\right)n\\
	\nonumber &=&\frac{\epsilon^2}{16}n\,.
	\end{eqnarray}
	Thus, for $\epsilon_1:=\frac{\epsilon^2}{64}$,
	$$\P\Big(2^{-1}\delta u-2\epsilon_1 n\leq L_{1+\delta}(n+C,n)-2n\Big)\leq \P\left(32^{-1}\epsilon^2 n\leq L_{1+\delta}(n+C,n)-\left((1+\delta)+\frac{1}{1+\delta}\right)n\right)\, .$$ 
	We have already seen that $L_{1+\delta}(0,n)-n/(1+\delta)$  has expectation $0$ and variance of order $n$, and also that $L_{1+\delta}(n+C,n)-L_{1+\delta}(0,n) - (1+\delta)n$ has expectation $C(1+\delta)$ and variance of order $n$, so we conclude that
	$$\lim_{n\to\infty}\P\left(32^{-1}\epsilon^2 n\leq L_{1+\delta}(n+C,n)-\left((1+\delta)+\frac{1}{1+\delta}\right)n\right)=0\,,$$
	and hence 
	$$\lim_{n\to\infty}\P\Big(2^{-1}\delta u-2\epsilon_1 n\leq L_{1+\delta}(n+C,n)-2n\Big)=0\,.$$
	To bound the second summand in \eqref{sum}, take $u=\epsilon n$ and write
	$$\lim_{n\to\infty} \P\Big( A(u)\leq 2^{-1}\delta u-\epsilon_1 n\Big) =
	\lim_{n\to\infty}  \P\Big( \frac{B( \epsilon n)}{n} + \frac{\min_{ z\geq \epsilon n } \bar B_{2^{-1}\delta}(z) }{n} \leq  -\epsilon_1 \Big) =0 .$$
	By \eqref{exit11}, this concludes the proof of 
	$$\lim_{n\to\infty}\P\left(Z_\nu(n+C,n)>\epsilon n \right)=0\,.$$
	To get the analog result for $\left\{Z_\nu(n+C,n)<-\epsilon n \right\}$ one  just needs to adapt the same argument.  
\begin{flushright}
$\blacksquare$.
\end{flushright}

\subsection{Local comparison and attractiveness}\label{sub:comparison}
In the next lemmas we will always construct $L_{\nu_1}$ and $L_{\nu_2}$ simultaneously with the basic coupling \eqref{basic}.

\begin{lemma}\label{comparison}
	If $x<y$ and $Z_{\nu_1}(y,n)\leq Z_{\nu_2}(x,n)$ then   
	$$L_{\nu_1}(y,n)-L_{\nu_1}(x,n)\leq L_{\nu_2}(y,n)-L_{\nu_2}(x,n)\,.$$
\end{lemma}

\textbf{Proof:}\\
	
	Recall the definition of the geodesic $\gamma(\xx,\yy)$ between two points $\xx < \yy$ in $\mathbb R \times \mathbb Z$ in Subsection \ref{sub:LPP}. 
	Denote by $\gamma_n^z(x)$ to the geodesic 
	between $(z,1)$ and $(x,n)$. Notice that
	$$L\left((z,1),(x,n)\right)=L\left((z,1),(y,m)\right)+L((y,m),(x,n))\,,$$
	for any $(y,m)\in\gamma_n^z(x)$. \\
	
	Assume that $Z_{\nu_1}(y,n)\leq Z_{\nu_2}(x,n)$, denote $z_1\equiv Z_{\nu_1}(y,)$ and $z_2\equiv Z_{\nu_2}(x,n)$. Let $\cc$ be a crossing point between the two geodesics $\gamma_n^{z_1}(y)$ and $\gamma_n^{z_2}(x)$. Such a crossing point always exists because $x\leq y$ and $z_1\leq z_2$ (by assumption). We remark that, by superaddivity of $L$,
	$$L_{\nu_2}(y,n)  \geq  \nu_2(z_2) + L\left((z_2,1),(y,n)\right) \geq  \nu_2(z_2) + L\left((z_2,1),\cc\right) + L\left(\cc,(y,n)\right)\,.$$
	We use this, and that (since $\cc\in\gamma_n^{z_2}(x)$)
	$$ \nu_2(z_2) + L\left((z_2,1),\cc\right)-L_{\nu_2}(x,n)= -L\left(\cc,(x,n)\right)\,,$$
	in the following inequality:
	\begin{eqnarray*}
		M^{(n)}_{\nu_2}(x,y)&=&L_{\nu_2}(y,n) - L_{\nu_2}(x,n) \\
		& \geq & \nu_2(z_2)+L\left((z_2,1),\cc\right) + L\left(\cc,(y,n)\right) - L_{\nu_2}(x,n)\\
		& = & L\left(\cc,(y,n)\right) - L\left(\cc,(x,n)\right)\,.
	\end{eqnarray*}
	By superaddivity,
	$$ - L\left(\cc,(x,n)\right)\geq L_{\nu_1}(\cc)-L_{\nu_1}(x,n)\,,$$
	and hence (since $\cc\in\gamma_{z_1}(y,n)$)
	\begin{eqnarray*}
		M^{(n)}_{\nu_2}(x,y)& \geq & L\left(\cc,(y,n)\right) - L\left(\cc,(x,n)\right)\\
		& \geq & L\left(\cc,(y,n)\right) + L_{\nu_1}(\cc)-L_{\nu_1}(x,n)\\
		& = & L_{\nu_1}(y,n)-L_{\nu_1}(x,n)\\
		& = & \Delta M^{(n)}_{\nu_1}(x,y)\,.
	\end{eqnarray*}
\begin{flushright}
	$\blacksquare$.
\end{flushright}

\begin{lemma}\label{lem:attract}
	Assume that $\nu_1(y)-\nu_1(x)\leq \nu_2(y)-\nu_2(x)$ for all $x<y$. Then 
	$$L_{\nu_1}(y,n)-L_{\nu_1}(x,n)\leq L_{\nu_2}(y,n)-L_{\nu_2}(x,n)\,,\,\forall\,x<y\,.$$
\end{lemma}

\textbf{Proof:}\\

	Denote 
	$$z_1:=Z_{\nu_1}(y,n)\,\,\mbox{ and }\,\,z_2:=Z_{\nu_2}(x,n)\,.$$
	If $z_1 \leq z_2$ then it follows from Lemma \ref{comparison} (we do not need to use the assumption). If  $z_1>z_2$ then
	\begin{eqnarray*}
		L_{\nu_2}(y,n)-L_{\nu_2}(x,n) -\big(L_{\nu_1}(y,n)-L_{\nu_1}(x,n)\big)&=&\\
		L_{\nu_2}(y,n)-\big(\nu_2(z_2)+L\left((z_2,1),(x,n)\right)\big)-\Big(\big(\nu_1(z_1)+L\left((z_1,1),(y,n)\right)\big)-L_{\nu_1}(x,n)\Big)&=&\\
		L_{\nu_2}(y,n)-\big(\nu_2(z_2)+L\left((z_1,1),(y,n)\right)\big)-\Big(\big(\nu_1(z_1)+L\left((z_2,1),(x,n)\right)\big)-L_{\nu_1}(x,n)\Big)&=&\\
		L_{\nu_2}(y,n)-\big(\nu_2(z_2)+L\left((z_1,1),(y,n)\right)\big)+\Big(L_{\nu_1}(x,n)- \big(\nu_1(z_1)+L\left((z_2,1),(x,n)\right)\big)\Big)&=&\\
		L_{\nu_2}(y,n)-\big(\nu_2(z_1)+L\left((z_1,1),(y,n)\right)\big)+\Big(L_{\nu_1}(x,n)- \big(\nu_1(z_2)+L\left((z_2,1),(x,n)\right)\big)\Big)&+&\\
		\big(\nu_2(z_1)-\nu_2(z_2)\big)-\big(\nu_1(z_1)-\nu_1(z_2)\big) \,.
	\end{eqnarray*}
	By super-additivity,
	$$L_{\nu_2}(y,n)-\big(\nu_2(z_1)+L_{z_1}(y,n)\big)\geq 0\, $$
	and 
	$$L_{\nu_1}(x,n)- \big(\nu_1(z_2)+L_{z_2}(x,n)\big)\geq 0\,,$$
	while, by assumption, 
	$$\nu_2(z_1)-\nu_2(z_2)\geq \nu_1(z_1)-\nu_1(z_2)\,,$$
	since $z_1>z_2$.
\begin{flushright}
	$\blacksquare$.
\end{flushright}	

\subsection{Proof of Theorem \ref{lppattract}}\label{sub:proof}
Without lost of generality we will assume that $\mu=1$ (again by Brownian scaling \eqref{scal}), and that $K=[0,C]$ with $C>0$. We take as an initial profile a Brownian motion with drift $1$,
$$B_1(x):=x+B(x)\,,$$
and also 
$$B_{\mu_\pm}:=\mu_{\pm}x+B(x)\,,$$
with  $\mu_\pm:=1\pm\delta$ and $\delta>0$. Thus,
$$B_{\mu_-}(y)-B_{\mu_-}(x)\leq B_{1}(y)-B_{1}(x)\leq B_{\mu_+}(y)- B_{\mu_+}(x)\,.$$
\begin{lemma}\label{control}
	Let $\mu\in(0,\infty)$ and assume \eqref{exitcontrol}. Then, for all $C>0$,
	$$\lim_{n\to\infty}\P\Big(Z_{\mu_-}(n+C,n)\leq Z_\nu(n,n)\,\mbox{ and }\,Z_{\nu}(n+C,n)\leq Z_{\mu_+}(n,n)\Big)=1\,.$$
\end{lemma}

\textbf{Proof:}\\

	Let us first prove that 
	$$\lim_{n\to\infty}\P\left(Z_{\mu_-}(n+C,n)\leq Z_\nu(n,n)\right)=1\,.$$ 
	For any $\epsilon>0$, 
	\begin{equation*}\P\left(Z_{\mu_-}(n+C,n) > Z_\nu(n,n)\right)\leq \P\left(Z_{\mu_-}(n+C,n) > -\epsilon n\right)+\P\left(Z_\nu(n,n)<-\epsilon n\right)\,.
	\end{equation*}
	Thus, by Lemma \ref{exit}, it is enough to show that (for fixed $\delta,C>0$) we can choose $\epsilon>0$ such that 
	\begin{equation}\label{eqcontrol}
	\lim_{n\to\infty}\P\left(Z_{\mu_-}(n+C,n)\leq -\epsilon n\right)=1\,.
	\end{equation}
	By shift invariance of Brownian Motion ($B(x+C)-B(C)\stackrel{dist.}{=}B(x)$), 
	$$\P\left(Z_{\mu_-}(n+C,n)>-\epsilon n\right)=\P\left( Z_{\mu_-}(n,n)>-\epsilon n-C\right)\leq\P\left(Z_{\mu_-}(n,n)>-2\epsilon n\right)\,,$$
	for $n\geq C/\epsilon$. Since 
	$$n=\mu^{-2}_{-}n+\left(1-\mu^{-2}_{-}\right)n\,\mbox{ and }\left(1-\mu^{-2}_{-}\right)<-\delta\,,$$
	(recall $\delta\in(0,1/2)$) by using shift invariance again,
	$$\P\left( Z_{\mu_-}(n,n)>-2\epsilon n\right) \leq \P\left(  Z_{\mu_-}(\mu^{-2}_{-}n,n) >(\delta-2\epsilon)n\right)\,.$$
	Hence, if $\epsilon<\delta/2$, Lemma \ref{exit} implies \eqref{eqcontrol}.  The proof of 
	$$\lim_{n\to\infty}\P\left( Z_{\nu}(n+C,n)\leq Z_{\mu_+}(n,n) \right)=1\,$$   
	is analogous.

\begin{flushright}
	$\blacksquare$.
\end{flushright}

If $Z_{\mu_-}(n+C,n)\leq Z_\nu(n,n)$ and $Z_{\nu}(n+C,n)\leq Z_{\mu_+}(n,n)$ then $Z_{\mu_-}(n+x,n)\leq Z_\nu(n,n)$ and $Z_{\nu}(n+x,n)\leq Z_{\mu_+}(n,n)$ for all $x\in[0,C]$. We use that $Z_\nu(y,n)$ is a non-decreasing function of $y$ (for fixed $n$). By Lemma \ref{comparison}, 
$$M_{\mu_-}^{(n)}(n,n+x)\leq M_\nu^{(n)}(n,n+x)\leq M_{\mu_+}^{(n)}(n,n+x)\,,$$
for all $x\in[0,C]$, and by Lemma \ref{lem:attract},
$$M_{\mu_-}^{(n)}(n,n+x)\leq M_1^{(n)}(n,n+x)\leq M_{\mu_+}^{(n)}(n,n+x)\,,$$
for all $x\in[0,C]$. Therefore,
\begin{eqnarray*}
	|M_\nu^{(n)}(n,n+x)-M_1^{(n)}(n,n+x)|&\leq & M_{\mu_+}^{(n)}(n,n+x)-M_{\mu_-}^{(n)}(n,n+x)\\
	&\leq & M_{\mu_+}^{(n)}(n,n+C)-M_{\mu_-}^{(n)}(n,n+C)\,,  
\end{eqnarray*}
for all $x\in[0,C]$. We use that $M_{\mu_+}^{(n)}(n,n+x)-M_{\mu_-}^{(n)}(n,n+x)$ is a non-decreasing function of $x$ (Lemma \ref{lem:attract}). Hence, if $Z_{\mu_-}(n+C,n)\leq Z^\nu(n,n)$ and $Z_{\nu}(n+C,n)\leq Z_{\mu_+}(n,n)$ then 
\begin{equation}\label{unifcontrol}
\sup_{x\in[0,C]}|M_\nu^{(n)}(n,n+x)-M_1^{(n)}(n,n+x)|\leq M_{\mu_+}^{(n)}(n,n+C)-M_{\mu_-}^{(n)}(n,n+C)\,.
\end{equation}

Since $M_{\mu_+}^{(n)}(n,n+C)-M_{\mu_-}^{(n)}(n,n+C)\geq 0$ (Lemma \ref{lem:attract}) and 
$$ \E\left(M_{\mu_+}^{(n)}(n,n+C)-M_{\mu_-}^{(n)}(n,n+C)\right)= (\mu_+-\mu_-)C=2\delta C\,,$$
we have that 
$$\P\left( M_{\mu_+}^{(n)}(n,n+C)-M_{\mu_-}^{(n)}(n,n+C)>\epsilon \right)\leq \frac{2C}{\epsilon}\delta\,.$$
Together with Lemma \ref{control} and \eqref{unifcontrol}, this implies that 
$$\limsup_{n\to\infty}\P\left( \sup_{x\in[0,C]}|M_\nu^{(n)}(n,n+x)-M_1^{(n)}(n,n+x)|>\eta \right)\leq \frac{2C}{\epsilon}\delta\,,$$
under \eqref{exitcontrol}. Since $\delta>0$ is arbitrary, we must have that 
$$\limsup_{n\to\infty}\P\left( \sup_{x\in[0,C]}|M_\nu^{(n)}(n,n+x)-M_1^{(n)}(n,n+x)|>\epsilon \right)=0 $$
under hypothesis \eqref{exitcontrol}, and Theorem \ref{lppattract} is proven. 
\begin{flushright}
$\blacksquare$
\end{flushright}

\section*{Conclusion}
 We proved that under mild conditions an initial flow passing through an infinite system of Brownian tandem queues converges in distribution to a Brownian Motion. The strong relationship between the queueing system and the Last Passage Brownian Percolation model is fundamental for the proof; since it allows to construct a coupling between different initial configurations using the concept of exit point in the LPP setting. This is a convenient way to manipulate the busy periods associated to the tandem queues. One wonders if this relation, or the one with the TABEP system, could be useful to compute non asymptotic formulae for the queueing system. An example of this kind of result, whose interpretation in the Brownian tandem queues setting has not yet been studied, is presented in \cite{FSW1}, where an explicit determinantal formula is obtained for the joint distribution of particles in a periodic finite system of particles interacting by one-sided reflection. \\

\textbf{Acknowledgements}
The authors would like to thank an anonymous referee for her helpful comments that greatly improved the presentation and clarity of this work.



\end{document}